\documentclass[10pt,a4paper]{article}
\usepackage[a4paper,top=2cm,bottom=2cm,left=2cm,right=2cm,marginparwidth=1.75cm]{geometry} 
\usepackage{graphicx}                      
\usepackage{placeins}
\usepackage{float} 
\usepackage{mathrsfs}
\usepackage{amsmath,amssymb}               
\usepackage{amsthm}
\usepackage{cite}                           
\usepackage{hyperref} 
\usepackage{caption}
\usepackage{subcaption}
\usepackage{titlesec}                       
\usepackage{float}                          
\usepackage{booktabs}                       
\usepackage{comment}
\newtheorem{Theorem}{Theorem}
\newtheorem{definition}{Definition}
\title{ Magnetic Resonance Dynamics via Fractional Bloch Equation: a Hybrid Computational Framework} 
\author{
    Neetu Garg$^{1}$, Diptiranjan Biswal$^{2}$, Varsha R$^{3}$ \\  
    {\small $^{1,2,3}$Department of Mathematics, National Institute of Technology Calicut, India} \\
    {\small Email: neetu@nitc.ac.in}
}
\date{}

\begin{document}

\maketitle

\begin{abstract}
Bloch equations are a powerful tool in describing the dynamics of nuclear magnetization in magnetic resonance phenomena. The fractional generalization of the Bloch equation effectively captures the anomalous relaxation and diffusion in porous, heterogeneous, and complex media. These equations describe how nuclear magnetization evolves under the influence of magnetic fields and relaxation processes. This work effectively employs a hybrid approach, the Laplace residual power series method, to investigate and analyze the fractional Bloch equation. A series solution is derived as the approximate solution for magnetization components. The influence of fractional order on each magnetization component in magnetization dynamics is analyzed and illustrated graphically. We conduct an error analysis to demonstrate the reliability and effectiveness of the proposed approach. The superiority of the suggested approach is shown using a comparative study with existing methods. The findings indicate the potential of the suggested approach as a reliable tool in understanding fractional magnetic resonance systems arising in applications such as NMR spectroscopy, MRI, MRF, and other complex heterogeneous materials. 
\end{abstract}
\section{Introduction}
Fractional calculus has emerged as a powerful mathematical framework that extends the classical calculus and deals with integrals and derivatives of arbitrary order \cite{podlubny1998fractional}. 
One of the key features of the fractional operators is their ability to incorporate the memory effect, in which the system's current state depends not only on current inputs but also on its entire past. This effect improves the modeling of real-life problems. This field is continuously evolving, and its relevance and applications in science and engineering are growing. Specifically, it has more significant applications in fields such as viscoelasticity, signal processing, biology, fluid dynamics, and control theory \cite{nonnenmacher1998applications,machado2011recent,miller1993introductio, diethelm2002analysis}.   
\par  In 1946, Felix Bloch introduced the Bloch equations, which are macroscopic equations that help us to understand the magnetic resonance as a set of first-order differential equations \cite{PhysRev.70.460}. Later, in 1956, Torrey added diffusion terms to these equations to study magnetization under the influence of shifting and relaxing magnetic fields \cite{torrey1956bloch}. The Bloch equations describe how nuclear magnetization evolves over time. 
They describe the spin system's relaxation dynamics by characterizing the time rate of change of magnetization. This system of equations has significant relevance in numerous scientific and engineering domains,  
particularly in applications such as nuclear magnetic resonance (NMR) \cite{singhi2021numerical}, magnetic resonance imaging (MRI) \cite{hazra2018numerical}, electron spin resonance (ESR) \cite{SINGH2017235}, and magnetic resonance fingerprinting (MRF) \cite{wang2019application}.  NMR imaging and spectroscopy technique is used to determine the structure of molecules. Scientists use it extensively in chemistry and biochemistry to understand how molecules are built and how they interact with each other. MRI is used to capture clear images of the interior of the human body. MRI is especially useful for examining the brain and detecting diseases such as cancer. ESR is used to study tiny particles called electrons in various materials. It also plays a role in developing quantum computers by enabling control over these electrons. Extremely accurate atomic clocks also use this technology to maintain precise timekeeping. This is crucial for applications like GPS, which rely on exact timing to provide accurate location information. MRF is a new and faster way to do MRI scans to get more detailed information about the tissues inside our body.
The  system of classical Bloch equations is given by \cite{kumar2014fractional}
\begin{equation}
\label{eq:bloch}
\left\{
\begin{aligned}
\frac{dM_x(t)}{dt} &= \omega_0 M_y(t) - \frac{M_x(t)}{T_2}, \\
\frac{dM_y(t)}{dt} &= -\omega_0 M_x(t) - \frac{M_y(t)}{T_2}, \\
\frac{dM_z(t)}{dt} &= \frac{M_0 - M_z(t)}{T_1}, \\
\end{aligned}
\right.
\end{equation}
with the initial state  
\[M_x(0) = 0, \quad M_y(0) = 100, \quad M_z(0) = 0.\]
Here, \(  M_x, M_y, ~\text{and }M_z \) are the components of the magnetization vector $M$ along the \( x \), \( y \), and \( z \) directions, respectively. 
\( M_0 \) represent the equilibrium magnetization along the \( z \)-axis. The angular frequency  $\omega_0$,
measured in radians, $\omega_0 = \gamma B_0$ and \(\omega_0 = 2\pi f_0\), where $B_0$ is the static magnetic field ($z$-component), $f_0$ is linear frequency measured in hertz,  and $\gamma$ is the gyromagnetic ratio. Then \(\frac{\gamma}{2\pi} = \frac{f_0}{B_0} = 42.57\) MHz/Tesla for water protons. Here, $T_1$ and $T_2$ represent the spin lattice and spin-spin relaxation times, respectively. The longitudinal relaxation time $T_1$ characterizes the process by which the nuclear spins regain the thermal equilibrium with the surrounding lattice through energy transfer. In contrast, $T_2$ describes the time over which the spins lose coherence in the transverse plane due to spin-spin interaction. 
The values of $T_1$ and $T_2$ vary across different body tissues.

The system of equations (\ref{eq:bloch}) admits the following analytical solution 
\begin{equation*}
\begin{aligned}
M_x(t) &= e^{-t/T_2} \left[ M_x(0)\cos(\omega_0 t) + M_y(0)\sin(\omega_0 t) \right], \\[6pt]
M_y(t) &= e^{-t/T_2} \left[ M_y(0)\cos(\omega_0 t) - M_x(0)\sin(\omega_0 t) \right], \\[6pt]
M_z(t) &= M_z(0)e^{-t/T_1} + M_0\left(1 - e^{-t/T_1} \right).
\end{aligned}
\end{equation*} 
Taking the asymptotic limit $t \rightarrow \infty$, 
the steady state solution associated with \eqref{eq:bloch} can be derived. 
Fractional Bloch equations are generalizations of the classical Bloch equations that capture memory effects, heterogeneity in the relaxation process, and complex magnetization dynamics. In this study, we investigate fractional Bloch equations described by \cite{kumar2014fractional}
 \begin{align}
    D_t^\alpha M_x (t) &= \omega_0 M_y (t) - \frac{M_x (t)}{T_2}, \label{eq:mx} \\
    D_t^\alpha M_y (t) &= -\omega_0 M_x (t) - \frac{M_y (t)}{T_2}, \label{eq:my} \\
    D_t^\alpha M_z (t) &= \frac{M_0 - M_z (t)}{T_1}, \label{eq:mz}
\end{align} with initial state \begin{equation} \label{IC}
    M_x(0) = 0 , ~M_y(0) = 100 ,~ \text{and} ~ M_z(0) = 0 .
\end{equation} 
Here $D_t^\alpha$ denotes the Caputo fractional derivative operator of order $\alpha\in (0,1]$ 
\cite{podlubny1998fractional}. 

In 2008, Magin et al. introduced this fractional Bloch equation to describe the magnetization dynamics in porous, complex, and heterogeneous media \cite{Magin09}.  This model enables the characterization of the anomalous relaxation and diffusion in such materials. Later, I. Petr\v{a}\v{s} studied the behavior and stability analysis of the fractional Bloch equation \cite{PETRAS2011341}. In 2011, Bhalekar introduced the time delay concept in the fractional Bloch equation with  \cite{bhalekar2011fractional}. 
Lately, many researchers have implemented mathematical methods for solving fractional Bloch equations, which include the operational matrix method \cite{SINGH2017235, mittal2019numerical},  homotopy perturbation method \cite{kumar2014fractional}, Lagrange's polynomial interpolation \cite{akgul2021new}, Laplace transform method \cite{ng}, Sumudu transform method \cite{singh2021new}, fractional variational iteration method \cite{prakash2019reliable}, analytical iterative method \cite{akshey2024approximate}, residual power series methods \cite{Sankar2025, DUBEY2022112691}, etc. 

The topic of solving fractional differential equations  has attracted significant attention, leading to the development of a variety of analytical and numerical methods. Recently, some authors have proposed hybrid methods that are developed by combining an integral transform with traditional analytical methods. One such combination is the Laplace transform and the residual power series method, known as the Laplace residual power series method (LRPSM) \cite{eriqat2020new}. This method delivers an accurate approximate solution as a convergent series. This method does not require perturbation, linearization, or discretization, unlike other existing methods. Recently, many scholars have employed the proposed technique to obtain solutions to a wide range of fractional differential equations, as shown in Refs. \cite{yadav2024constructing, bekhouche2025coupling,Yadav2025,garg2026}. 

This paper explains the implementation of LRPSM to find approximate series solutions of the fractional Bloch equations. We convert the system of fractional differential equations in to an equivalent system of algebraic equations and determine the unknown coefficients in the series solution in the Laplace domain. We demonstrate a rigorous analysis of the obtained approximate solution of magnetization components through graphical illustrations. We discuss the phase-plane trajectories describing the transverse components of magnetization and the influence of fractional orders on the Bloch equations. Furthermore, to highlight the superiority of the proposed method, we present a comparative study with other existing methods. 

\par This remainder of this paper is organized in the following manner:  Section \ref{sec2} presents preliminary concepts related to Caputo fractional derivative and the suggested method.  In Section \ref{sec3}, we rigourously outlines the methodology of the suggested approach for deriving the solution of a system of equations. The detailed solution process for the fractional Bloch equation using LRPSM is discussed in Section \ref{sec4}.  We analyze the numerical findings and compare them with other existing methods in Section \ref{sec5}. Finally, we conclude in Section \ref{sec6}.
\section{Mathematical Preliminaries} \label{sec2}
This section is devoted to the fundamental definitions and features of the Caputo fractional derivative and the Laplace transform. 
Additionally, we recall theorems related to fractional power series. 
\begin{definition}[Caputo derivative]
The Caputo fractional derivative with order $\alpha > 0$ of the given
function $\varphi(t), ~t \in (a,b)$ and $m-1 < \alpha \leq m ~ (m\in \mathbb{Z}^+)$, is defined as \cite{podlubny1998fractional}
\[ D_t^\alpha \varphi(t)=
 \frac{1}{\Gamma(m - \alpha)} \int_a^t (t - s)^{m-\alpha-1} \varphi^{(m)}(s) \, ds.
\]
\end{definition}
\begin{definition}[Laplace transform]
    The Laplace transform of a piecewise continuous function $\varphi(t)$ of exponential order $\delta$ is defined as follows:
    \begin{equation*}
        \Phi(s)=\mathcal{L}\{\varphi(t)\} = \int _0 ^\infty e^{-st} \varphi(t) \, dt, ~~ s>\delta.
    \end{equation*}
    The inverse Laplace transform of the function $\Phi(s)$ is defined as,
    \begin{equation*}
        \varphi(t)= \mathcal{L}^{-1}\{\Phi(s) \} = \int _{\eta-i\infty} ^{\eta+i\infty} e^{st} \Phi(s) \, ds, ~~~ \eta =Re(s)>0.
    \end{equation*}
\end{definition}
\noindent
Some characteristic features of the Laplace transform are listed below:\cite{eriqat2020new}
\begin{enumerate}
    \item \( \lim\limits_{s \to \infty} s \Phi(s) = \lim\limits_{t \to 0}\varphi(t) = \varphi(0). \)
    \item \(\mathcal{L}\{D_t^\alpha \varphi(t)\} = s^\alpha \Phi(s) - \sum\limits_{k=0}^{m-1} s^{\alpha - k - 1} \varphi^{(k)}(0), \quad m - 1 < \alpha \leq m.\)
\end{enumerate}

\begin{Theorem}
\cite{eriqat2020new} Let  $ {\varphi} (t) $ admits a fractional power series (FPS) representation about \( t = 0 \) given by
\begin{equation} \label{fps}
{\varphi}(t) = \sum_{n=0}^{\infty} c _n t^{n\alpha}, \quad 0 < m - 1 < \alpha \leq {m}, \quad 0 \leq t < b.
\end{equation}
Suppose \( \varphi(t)  \) is continuous and   \( D_{t}^{n\alpha} \varphi(t) 
\) exists and continuous on  $t\in (0,b)$ for \( n \geq 0\). 
Then the coefficients \( c_n \) in Eq. \eqref{fps} take the form
\[
c_n = \frac{D_{t}^{n\alpha} \varphi(0)}{\Gamma(n\alpha + 1)}, \ n =0,1,2, \hdots  
\]

\end{Theorem}
\begin{Theorem}
\cite{eriqat2020new} Assume that the Laplace transform $\Phi(s) =\mathcal{L}\{\varphi(t)\}$ admits a FPS representation of the form
\[
\Phi(s) = \sum_{n=0}^{\infty} \frac{c_n}{s^{n\alpha+1}}, \quad 0 < \alpha \leq 1, \ s > 0,
\]
where $c_n = D_{t}^{n\alpha} \varphi (0)$. 
 The inverse Laplace transform of the FPS of $\Phi(s)$ has the following form
\[
\varphi(t) = \sum_{n=0}^{\infty} \frac{D_{t}^{n\alpha} \varphi(0)}{\Gamma(n\alpha + 1)} t^{n\alpha}, \quad 0 < \alpha \leq 1, \ t \geq 0.\]

\end{Theorem}
\begin{Theorem}
    \cite{eriqat2020new} Suppose $\Phi(s)= \mathcal{L}\{\varphi(t) \} $ have the FPS representation in the form $\Phi(s)= \sum\limits_{n=0}^\infty \frac{c_n}{s^{n\alpha +1}}, ~ 0 < \alpha \leq 1,$ and $\left| s\mathcal{L}\{ D_t^{(n+1)\alpha} \varphi(t) \} \right| \leq \mathrm{M} ~~\left( \mathrm{M} \in \mathbb{R}^+ \right)$ on $0 \leq s \leq d$. Then, the remainder $\mathcal{R}_n(s)$ holds the following bound 
 \begin{equation*}
    \left| \mathcal{R}_n(s) \right| \leq \frac{\mathrm{M}}{s^{(n+1)\alpha + 1}}, \quad 0 \leq s \leq d, \quad \left( d \in \mathbb{R}^+ \right).
\end{equation*}
\end{Theorem}
\noindent
\section{Methodology} \label{sec3}
This section briefly presents a concise overview of the LRPSM approach for a system of equations, which will be employed to derive an approximate solution to the fractional Bloch equation. For a detailed methodology, we consider a general system of $\ell$ fractional differential equations involving the Caputo derivatives of order $\alpha\in (0,1]$ 
\begin{equation} \label{eq1}
D_t^{\alpha} \varphi_i(t) = f_i( t, \varphi_1(t),\varphi_2(t), \hdots \varphi_\ell (t)), ~ t\geq 0, ~ i=1, 2, \hdots \ell,\end{equation}
subject to the initial condition  $\varphi_i(0) = \varphi_{i,0}$.
~~Here, $f_i$'s are known continuous real-valued functions and $\varphi_i$ denotes the unknown smooth functions to be determined. For each $i=1,2, \hdots \ell$, we follow the below stated steps.

\noindent Step 1: Perform the Laplace transform on Eqs. \ref{eq1} and simplify using the initial condition. Then, we obtain
\begin{equation} \label{eq7}
    \Phi_i(s) = \frac{\varphi_{i,0}}{s} + \frac{1}{s^\alpha}\mathcal{L}\left\{f_i \left(t,\varphi_1(t)), \varphi_2(t), \hdots \varphi_\ell(t)\right) \right\},
\end{equation} where $\Phi_i(s) = \mathcal{L}\{\varphi_i(t)\}$. \\
Step 2: Suppose that the solutions $\Phi_i(s)$ of the Eq. \eqref{eq7} admit the following series expansion
    \begin{equation}
    \Phi_i(s) = \sum_{n=0}^{\infty} \frac{c_{i,n}}{s^{n\alpha + 1}}, \quad s > 0, \quad i=1,2, \hdots \ell.
    \end{equation}
    The corresponding $k$-th truncated series, denoted by $\Phi_{i}^k$, is given by, \begin{equation}
        \Phi_i^k(s) = \sum_{n=0}^{k} \frac{c_{i,n}}{s^{n\alpha + 1}}= \frac{\varphi_{i,0}}{s}+ \sum_{n=1}^{k} \frac{c_{i,n}}{s^{n\alpha + 1}}, \quad s > 0,\quad i=1,2, \hdots \ell.
    \end{equation}
     
\noindent Step 3: Define the Laplace residual function $\mathcal{L}\text{Res}_{\Phi_i}(s)$ associated with $\Phi_i(s)$ and consider its  $k$-th Laplace residual function $\mathcal{L}\text{Res}_{\Phi_i}^k(s)$ as follows
\begin{equation*}
\mathcal{L}\text{Res}_{\Phi_i}(s) = \Phi_i(s) - \frac{\varphi_{i,0}}{s} - \frac{1}{s^\alpha} \mathcal{L}\left\{ f_i\left(t, \mathcal{L}^{-1}\{ \Phi_1(s)\}, \mathcal{L}^{-1}\{ \Phi_2(s)\}, \hdots \mathcal{L}^{-1}\{ \Phi_\ell(s) \}\right) \right\},
\end{equation*}
\begin{equation} \label{eq10}
\mathcal{L}\text{Res}_{\Phi_i}^k(s) = \Phi_i^k(s) - \frac{\varphi_{i,0}}{s} - \frac{1}{s^\alpha} \mathcal{L}\left\{ f_i\left(t, \mathcal{L}^{-1}\{ \Phi_1^k(s)\}, \mathcal{L}^{-1}\{ \Phi_2^k(s)\}, \hdots \mathcal{L}^{-1}\{ \Phi_\ell^k(s) \}\right) \right\}.
\end{equation}
Step 4: Multiply each side of Eq. \eqref{eq10} by $s^{k\alpha +1}$. \\
Step 5: Find all $c_{i,n}$ (for $n=1, 2, 3, \dots$) by using
    \[
    \lim_{s \to \infty} s^{k\alpha + 1} \mathcal{L}\text{Res}_{\Phi_i} ^k(s) =0, \quad k = 1, 2, 3, \dots
    \]
Step 6: Determine $\Phi_i^k(s)$ by substituting the values of $c_{i,n}$ into the $k$-th approximate solution.\\
Step 7: Find the $k$-th approximate solution $\varphi_i^k(t)$ by taking the inverse Laplace transform of $\Phi_i^k(s)$.

\section{LRPSM Implementation on the Fractional Bloch Equations}\label{sec4}
This section describes the implementation of the LRPSM on the fractional Bloch equations. Consider the system of equations \eqref{eq:mx}- \eqref{IC}.
~~Suppose \begin{equation*}
    \mathcal{L}\{M_x(t)\} = \mathscr{M}_x(s), ~ \mathcal{L}\{M_y(t)\} = \mathscr{M}_y(s), ~\mathcal{L}\{M_z(t)\} = \mathscr{M}_z(s). ~ 
\end{equation*}
Apply Laplace transform to the system \eqref{eq:mx}-\eqref{eq:mz} and 
employing the initial condition, we obtain
\begin{equation}
    \begin{split}
\mathscr{M}_x(s) =& ~\frac{\omega_0}{s^\alpha} \mathscr{M}_y(s) - \frac{1}{T_2 s^\alpha} \mathscr{M}_x(s),\\
 \mathscr{M}_y(s) =& ~\frac{100}{s} -\frac{\omega_0}{s^\alpha} \mathscr{M}_x(s) - \frac{1}{T_2 s^\alpha} \mathscr{M}_y(s),\\
 \mathscr{M}_z(s) = & ~\frac{M_0}{ T_1 s^{\alpha+1}}  - \frac{\mathscr{M}_z(s)}{T_1 s^\alpha}.
    \end{split}
\end{equation}
According to LRPSM, we assume that the functions $\mathscr{M}_x(s), ~ \mathscr{M}_y(s), ~ \text{and} ~\mathscr{M}_z(s)$  have the following FPS representations
\begin{equation}
\label{eq13}
    \begin{split}
        \mathscr{M}_x(s)=&  
 \sum_{n=0}^{\infty} \frac{\mathfrak{a}_n}{s^{n\alpha+1}}, \quad  
 \mathscr{M}_x^k(s) = \sum_{n=0}^{k} \frac{\mathfrak a_n}{s^{n\alpha+1}},\\
\mathscr{M}_y(s)= &\sum_{n=0}^{\infty} \frac{\mathfrak{b}_n}{s^{n\alpha+1}}, \quad  
 \mathscr{M}_y^k(s) = \sum_{n=0}^{k} \frac{\mathfrak b_n}{s^{n\alpha+1}},\\
\mathscr{M}_z(s)=& \sum_{n=0}^{\infty} \frac{\mathfrak{c}_n}{s^{n\alpha+1}}, \quad  
 \mathscr{M}_z^k(s) = \sum_{n=0}^{k} \frac{\mathfrak c_n}{s^{n\alpha+1}},
    \end{split}
\end{equation} where $ \mathscr{M}_x^k,  ~\mathscr{M}_y^k,$ and $ \mathscr{M}_z^k(s)$ are the $k$-th truncated series of $ \mathscr{M}_x(s), ~  \mathscr{M}_y(s),$ and $ \mathscr{M}_z(s)$, respectively.
Here,
\begin{equation}
    \mathfrak a_0 = \lim_{s \to \infty} s\mathscr{M}_x(s) = M_x(0)=0, ~ \mathfrak b_0 = \lim_{s \to \infty} s \mathscr{M}_y(s) = M_y(0)=100,~\text{ and } ~ \mathfrak c_0 = \lim_{s \to \infty} s \mathscr{M}_z(s) = M_z(0)=0,
\end{equation}
Now, we define the Laplace residual functions of $\mathscr{M}_x(s), ~ \mathscr{M}_y(s), ~ \text{and} ~\mathscr{M}_z(s)$ as follows
\begin{equation}
\begin{split}
\mathcal{L}{\text{Res}_{\mathscr{M}_x}}(s) =&~ \mathscr{M}_x(s) - \frac{\omega_0}{s^\alpha} \mathscr{M}_y(s) + \frac{1}{T_2 s^\alpha} \mathscr{M}_x(s),\\
\mathcal{L}{\text{Res}_{\mathscr{M}_y}}(s) = &~\mathscr{M}_y(s) - \frac{\mathfrak b_0}{s} + \frac{\omega_0}{s^\alpha} \mathscr{M}_x(s) + \frac{1}{T_2 s^\alpha} \mathscr{M}_y(s),\\
    \mathcal{L}\text{Res}_{\mathscr{M}_z}(s) =&~ \mathscr{M}_z(s) - \frac{M_0}{T_1 s^{\alpha+1}}  + \frac{1}{T_1 s^\alpha}\mathscr{M}_z(s).
\end{split}
\end{equation}
Thus, the corresponding k-th residual functions are given as \begin{equation} \label{eq17}
    \begin{split}
    \mathcal{L}{\text{Res}^{k}_{\mathscr{M}_x}}(s) = &~{\mathscr{M}_x^k}(s) - \frac{\omega_0}{s^\alpha} {\mathscr{M}_y^k}(s) + \frac{1}{T_2 s^\alpha} {\mathscr{M}_x^k}(s),\\
    \mathcal{L}{\text{Res}^{k}_{\mathscr{M}_y}}(s) =&~ {\mathscr{M}_y^k}(s) - \frac{\mathfrak b_0}{s} + \frac{\omega_0}{s^\alpha} {\mathscr{M}_x^k}(s) + \frac{1}{T_2 s^\alpha} {\mathscr{M}_y^k}(s),\\
        \mathcal{L}\text{Res}_{\mathscr{M}_z}^{k}(s) =& ~\mathscr{M}_z^k(s) - \frac{M_0}{T_1 s^{\alpha+1}}  +\frac{1}{T_1 s^\alpha} \mathscr{M}_z^k(s).
    \end{split}
\end{equation}

To compute the value of $\mathfrak a_1,~ \mathfrak b_1$, and $\mathfrak c_1$, we consider $k=1$ in Eq. \eqref{eq17} and utilize \eqref{eq13} to obtain
\begin{equation} \label{eq18}
    \begin{split}
        \mathcal{L}{\text{Res}^{1}_{\mathscr{M}_x}}(s) = ~&~\left(\frac{\mathfrak a_0}{s} + \frac{\mathfrak a_1}{s^{\alpha+1}} \right)- \frac{\omega_0}{s^\alpha} \left( \frac{\mathfrak b_0}{s} + \frac{\mathfrak b_1}{s^{\alpha+1}} \right) + \frac{1}{T_2 s^\alpha} \left( \frac{\mathfrak a_0}{s} + \frac{\mathfrak a_1}{s^{\alpha+1}} \right),\\ 
        \mathcal{L}{\text{Res}^{1}_{\mathscr{M}_y}}(s) =&~ \left( \frac{\mathfrak b_0}{s} + \frac{\mathfrak b_1}{s^{\alpha+1}}\right) - \frac{\mathfrak b_0}{s} + \frac{\omega_0}{s^\alpha} \left( \frac{ \mathfrak a_0}{s} + \frac{\mathfrak  a_1}{s^{\alpha+1}} \right) + \frac{1}{T_2 s^\alpha} \left( \frac{\mathfrak b_0}{s} + \frac{\mathfrak b_1}{s^{\alpha+1}} \right),\\
        \mathcal{L}\text{Res}_{\mathscr{M}_z}^{1}(s) &= ~ 
\left( \frac{\mathfrak c_0}{s} + \frac{\mathfrak c_1}{s^{\alpha+1}} \right) - \frac{M_0}{T_1 s^{\alpha+1}}  + \frac{1}{T_1 s^\alpha} \left( \frac{\mathfrak c_0}{s} + \frac{\mathfrak c_1}{s^{\alpha+1}} \right).
    \end{split}
\end{equation} Multiplying each equation of the system \eqref{eq18} with $s^{\alpha+1}$ and taking the limit $s \rightarrow \infty$, yields
\begin{equation}
    \begin{split}
    \lim_{s \to \infty}s^{\alpha+1}\mathcal{L}\text{Res}_{\mathscr{M}_x}^{1}(s) &= \mathfrak a_1 -\omega_0 \mathfrak b_0=0,\\
    \lim_{s \to \infty}s^{\alpha+1}\mathcal{L}\text{Res}_{\mathscr{M}_y}^{1}(s) &= \mathfrak b_1 + \frac{\mathfrak b_0}{T_2} =0,\\
        \lim_{s \to \infty}s^{\alpha+1}\mathcal{L}\text{Res}_{\mathscr{M}_z}^{1}(s) &= \mathfrak c_1- \frac{M_0}{T_1}= 0.
    \end{split}
\end{equation}
Hence, we derive $\mathfrak a _1 = \omega _0 \mathfrak b_0, ~ \mathfrak b_1 = - \frac{\mathfrak b_0}{T_2}$, and $\mathfrak c_1 = \frac{M_0}{T_1}$.
Thus, continuing in this manner for $k=2,3, 4, \hdots$ we derive 
\[\mathfrak a _2 = -\frac{2 \omega_0 \mathfrak b_0}{T_2},\ \mathfrak a_3 = -\omega_0^3 \mathfrak b_0 + \frac{3 \omega_0 \mathfrak b_0}{T_2^2}, \ \mathfrak a_4 = \frac{4 \omega_0^3 \mathfrak b_0}{T_2} - \frac{4 \omega_0 \mathfrak b_0}{T_2^3}, ~ \hdots\]
\[\mathfrak b_2 =-\omega_0^2 \mathfrak b_0 + \frac{\mathfrak b_0}{T_2^2},\ \mathfrak b_3 = \frac{3\omega_0^2 \mathfrak b_0}{T_2} - \frac{\mathfrak b_0}{T_2^3},\ \mathfrak b_4 = \omega_0^4 \mathfrak b_0 - \frac{6 \omega_0^2 \mathfrak b_0}{T_2^2} + \frac{\mathfrak b_0}{T_2^4},~ \hdots\]
and 
\[\mathfrak c_2 = -\frac{M_0}{T_1^2},\
\mathfrak c_3  = \frac{M_0}{T_1^3},\
\mathfrak c_4 = -\frac{M_0}{T_1^4}, ~ \hdots\]
 By taking the inverse Laplace transform of $\mathscr{M}_x, ~ \mathscr{M}_y$, and $\mathscr{M}_z,$ we obtain approximate series solutions in the form
\[\mathring{M}_x(t)
 = 100 \, \omega_0 \left( \frac{t^\alpha}{\Gamma(\alpha+1)} - \frac{2t^{2\alpha}}{T_2 \Gamma(2\alpha+1)} + \frac{3t^{3\alpha}}{T_2^2 \Gamma(3\alpha+1)} - \frac{\omega_0^2 t^{3\alpha}}{\Gamma(3\alpha+1)} + \dots \right)
,\]
\begin{equation*}
    \begin{split}
\mathring{M}_y(t)= 100 E_\alpha \left( -\frac{t^\alpha}{T_2}\right) - \frac{100 \omega_0^2 t^{2\alpha}}{\Gamma(2\alpha+1)} + \frac{100 \omega_0^2}{T_2} \left( \frac{2t^{2\alpha}}{\Gamma(3\alpha+1)} + \frac{t^{3\alpha}}{\Gamma(3\alpha+1)} \right) + \dots,
    \end{split}
    \end{equation*}
and 
\[ 
\mathring{M}_z 
=  M_0 \left( 1 - E_{\alpha} \left( -\frac{t^\alpha}{T_1} \right) \right),\] where $E_\alpha(t)$ denotes the Mittag-Leffler function of one parameter. Here, $\mathring{M}_x(t),~ \mathring{M}_y(t),$ and $\mathring{M}_z(t)$ denote the  approximate solutions of magnetization components $M_x(t),~ M_y(t),$ and $M_z(t),$ respectively.

\section{Results and Discussions}\label{sec5}
This section presents numerical results obtained using the proposed method and discusses the physical dynamics of the approximate solution. To assess  the accuracy and reliability of  
the proposed method, two types of error measures are considered, as described below
\begin{itemize}
    \item[] Absolute errors $\mathbf{e}_x(t) =  \left| M_x(t) - \mathring{M}_x(t) \right|$, $\mathbf{e}_y(t) = \left| M_y(t) - \mathring{M}_y(t) \right|$, $\mathbf{e}_z(t) = \left| M_z(t) - \mathring{M}_z(t) \right|$, \\
    \item[] Relative error $= \frac{\mathbf{e}_i(t)}{M_i(t)}, $ where $i \in \{x,y,z \}$.
\end{itemize}
We assume $\omega_0  = 1, ~T_1=1~(ms)^q$, and $T_2=20~(ms)^q$. 
In Table \ref{tab:comparison},  we consider the approximate solution by taking 10 terms for $\alpha =1$. We observe that the approximate solution is 
consistent with the solution obtained by other methods. The exact and approximate values are closely same. We observe that this method is more efficient than other existing methods, which indicate the reliability of the suggested method. 
\begin{table}[h]
\centering
\caption{Comparison with existing methods.} 
\resizebox{\textwidth}{!}{%
\begin{tabular}{ccccccccccc}
\toprule
{} & $t$ & Exact solution & LRPSM &ARA-RPSM \cite{Sankar2025}& FVIM\cite{prakash2019reliable} & FHPTM \cite{prakash2019reliable}& \begin{tabular}[c]{@{}c@{}}Operational\\ Matrix \cite{SINGH2017235}\end{tabular} & \begin{tabular}[c]{@{}c@{}}HPM\cite{kumar2014fractional}\end{tabular} &  \begin{tabular}[c]{@{}c@{}}ATIM\cite{akshey2024approximate}\end{tabular} \\
\midrule
$M_x(t)$& 0.1 & 9.9335 & 9.9335 & 9.9335  & 9.9335 & 9.95256 & 9.9245 & 9.9335  & 9.9335 \\
& 0.3 & 29.1120 & 29.1120 & 29.1120 & 29.1120 & 29.6219 & 29.1080 & 29.1034  & 29.1034 \\
& 0.5 & 46.7589 & 46.7588 & 46.7588 & 46.7589 & 49.0948 & 46.7732 & 46.6823  & 46.6823 \\
& 0.7 & 62.2060 & 62.2060  & 62.2060 & 62.2074 & 68.5244 & 62.2180 & 61.8762 &  61.8762 \\
& 0.9 & 74.8859 & 74.8859 & 74.8859 & 74.8943 & 88.0717 & 74.8814 & 73.8911  & 73.8911 \\
\midrule
$M_y(t)$& 0.1 & 99.0042 & 99.0042 & 99.0042 & 99.0042 & 99.0037 & 99.0213 & 99.0187  &99.0030 \\
& 0.3 & 94.1113 & 94.1113 & 94.1113 & 94.1113 & 94.0731 & 94.1645 & 94.1837 &  94.0573 \\
& 0.5 & 85.5915 & 85.5915& 85.5915 & 85.5914 & 85.2988 & 85.5689 & 85.5518 &  85.2445 \\
& 0.7 & 73.8536 & 73.8536 & 73.8536 & 73.8529 & 72.7400 & 73.7886 & 73.1630 &  72.6465 \\
& 0.9 & 59.4258 & 59.4258  & 59.4258 & 59.4215 & 56.4183 & 59.3742 & 57.0572 &  56.3451 \\
\midrule
$M_z(t)$& 0.1 & 0.0952 & 0.0952 & 0.0952 & 0.0952 & 0.0952 & 0.0952 & 0.0952 & 0.0952 \\
& 0.3 & 0.2592 & 0.2592 & 0.2592 & 0.2592 & 0.2592 & 0.2592 & 0.2592 &  0.2592 \\
& 0.5 & 0.3935 & 0.3935 & 0.3935 & 0.3935 & 0.3935 & 0.3935 & 0.3935 &  0.3935 \\
& 0.7 & 0.5034 & 0.5034 & 0.5034 & 0.5034 & 0.5034 & 0.5034 & 0.5034 &  0.5034 \\
& 0.9 & 0.5934 & 0.5934 & 0.5934 & 0.5934 & 0.5934 & 0.5934 & 0.5934 &  0.5934 \\
\bottomrule
\end{tabular}%
}
\label{tab:comparison}
\end{table}

\begin{table}[h!]
\centering
\caption{Absolute and relative errors. }
\begin{tabular}{cccc}
\toprule
\textbf{} & {\( t \)} & {Absolute Error} & {Relative Error} \\
\midrule
 $M_x(t)$& 0.1 & 0 & 0 \\
& 0.2 & 5.32907051820075e-14 & 2.70934046325674e-15 \\
& 0.3 & 3.76942921320733e-12 & 1.29480040762691e-13 \\
& 0.4 & 8.87325768417213e-11 & 2.32462326405266e-12 \\
& 0.5 & 1.02684083458371e-09 & 2.19603535751521e-11 \\
\midrule
 $M_y(t)$& 0.1 & 1.42108547152020e-14 & 1.43537961139099e-16 \\
& 0.2 & 5.68434188608080e-14 & 5.85824534945981e-16 \\
& 0.3 & 2.41584530158434e-12 & 2.56700770424182e-14 \\
& 0.4 & 5.85345105719171e-11 & 6.48349962837456e-13 \\
& 0.5 & 6.89183821123152e-10 & 8.05201269150381e-12 \\
\midrule
 $M_z(t)$& 0.1 & 5.55111512312578e-17 & 5.83329603774665e-16 \\
& 0.2 & 5.55111512312578e-16 & 3.06235901422036e-15 \\
& 0.3 & 4.34097202628436e-14 & 1.67487546296745e-13 \\
& 0.4 & 1.01674224595172e-12 & 3.08402811188705e-12 \\
& 0.5 & 1.17418297307381e-11 & 2.98417907788256e-11 \\
\bottomrule
\end{tabular}
\label{tab:errors}
\end{table}

\noindent The absolute and relative errors by taking 10 terms of series solution for $\alpha =1$ are summarized in Table \ref{tab:errors}. We observe that both errors remain negligibly small throughout the considered interval. The absolute errors increase over time. The consistently small relative errors indicate that the method provides stable, reliable approximations for the magnetization components.

\begin{figure}[H]
\begin{subfigure}[H]{0.5\textwidth}
    \centering
   \includegraphics[width=\linewidth]{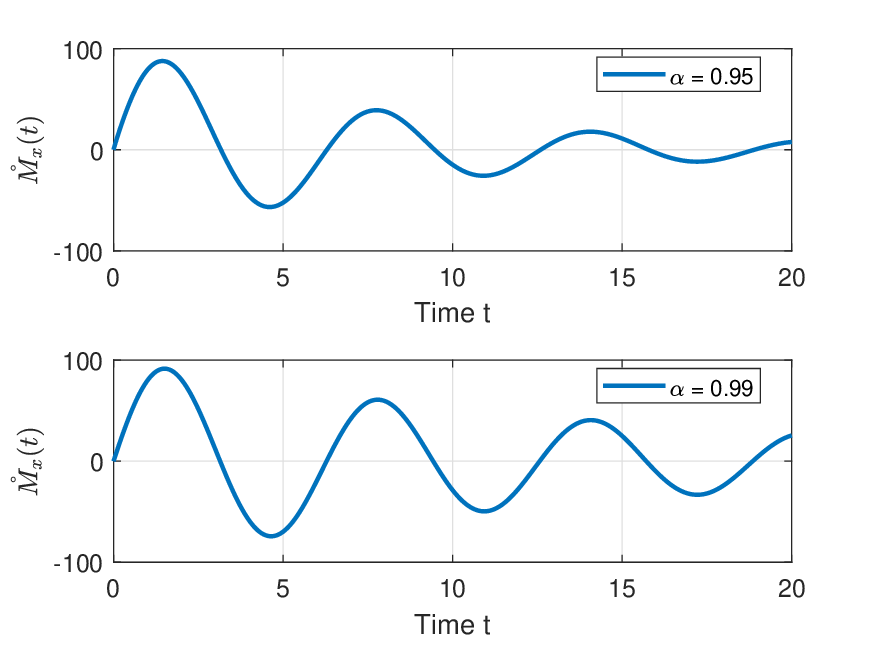}
    \caption{$\mathring{M}_x(t)$ }
    \label{fig:M_xdiffalpha}
     \end{subfigure}
    \hfill
\begin{subfigure}[H]{0.5\textwidth}
        \centering
    \includegraphics[width=\linewidth]{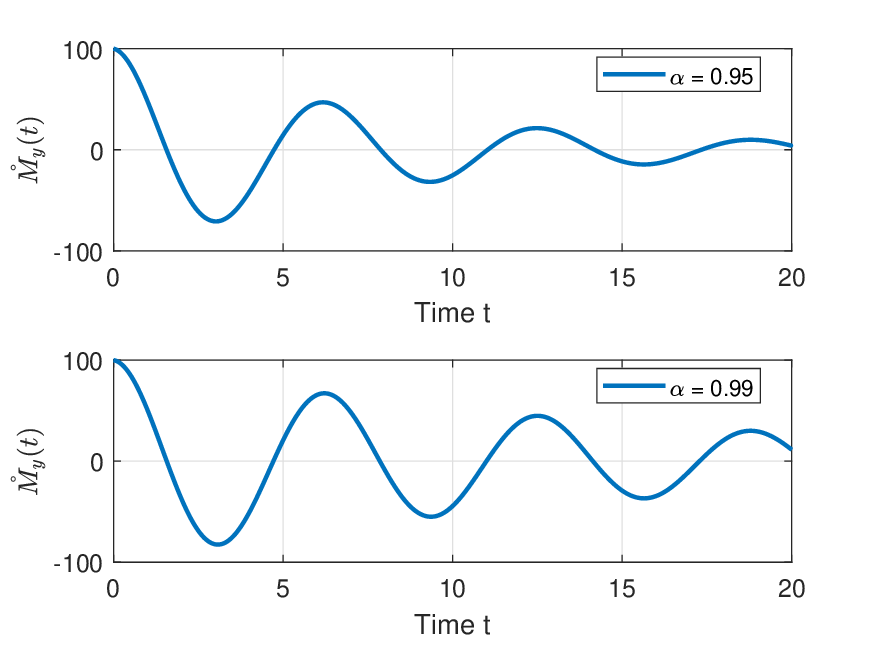}
    \caption{$\mathring{M}_y(t)$}
    \label{fig:M_ydiffalpha}
    \end{subfigure}
    \caption{Time evolution of the transverse magnetization components $\mathring{M}_x(t)$ and $\mathring{M}_y(t)$ for different fractional orders (with 30 terms in series solution).}
    \label{fig:M_xydiffalpha}
\end{figure}

The influence of fractional orders on transverse magnetization components $M_x(t)$ and $M_y(t)$ is depicted in Figure  \ref{fig:M_xydiffalpha}. We observe that both  transverse components exhibit  damped oscillatory behavior over time, which is a characteristic feature of transverse relaxation in the classical Bloch system. Moreover, the degree of damping is strongly influenced by the fractional order. As $\alpha$ decreases from 1, the damping increases and leads to a faster attenuation,  
particularly as time increases. 
The variation in decay behavior according to the fractional order highlights that the fractional Bloch system  effectively delivers the memory-dependent relaxation effects. 

\begin{figure}[h]
    \centering
    \includegraphics[width=0.5\linewidth]{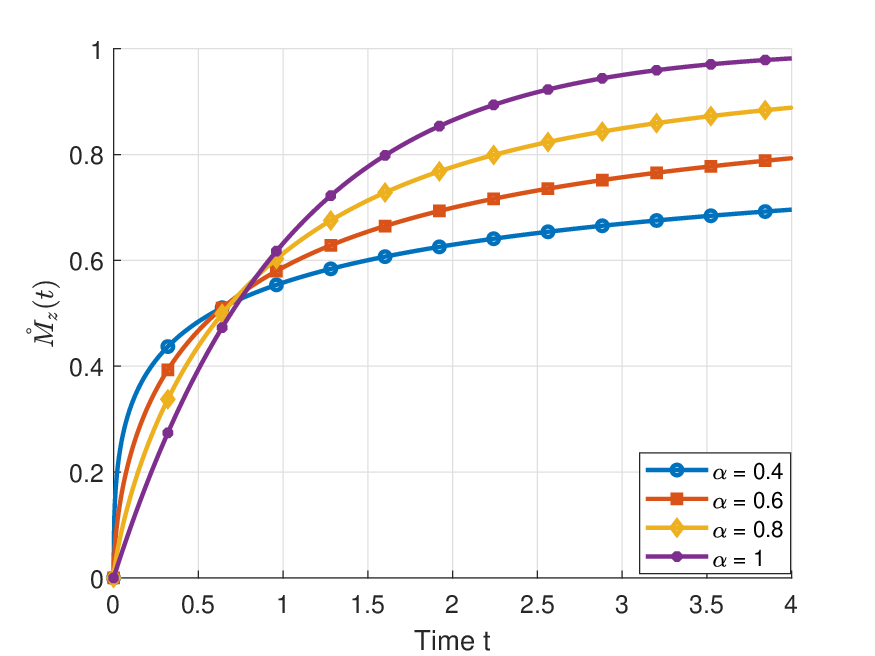}
    \caption{Time evolution of longitudinal magnetization component $\mathring{M}_z(t)$ for different fractional orders (with 30 terms in the series solution).}
    \label{fig:M_zdiffalphas}
\end{figure}

 \noindent The variation of longitudinal magnetization component $M_z(t)$ with respect to time for various fractional orders in the relaxation dynamics is demonstrated in Figure \ref{fig:M_zdiffalphas}. The magnetization increases monotonically and approaches to a steady state, indicating spin-lattice relaxation. A notable feature of the component $M_z(t)$ is observed  near time $t\approx
1$, where all curves corresponding to different fractional orders intersects.  Initially, we observe that $M_z(t)$ increases as $\alpha$ decreases. However, beyond this certain time, this trend reverse and $M_z(t)$ decreases with decrease in  $\alpha$. 

\begin{figure}[H]
\centering

\begin{subfigure}[b]{0.48\textwidth}
    \centering
    \includegraphics[width=\linewidth]{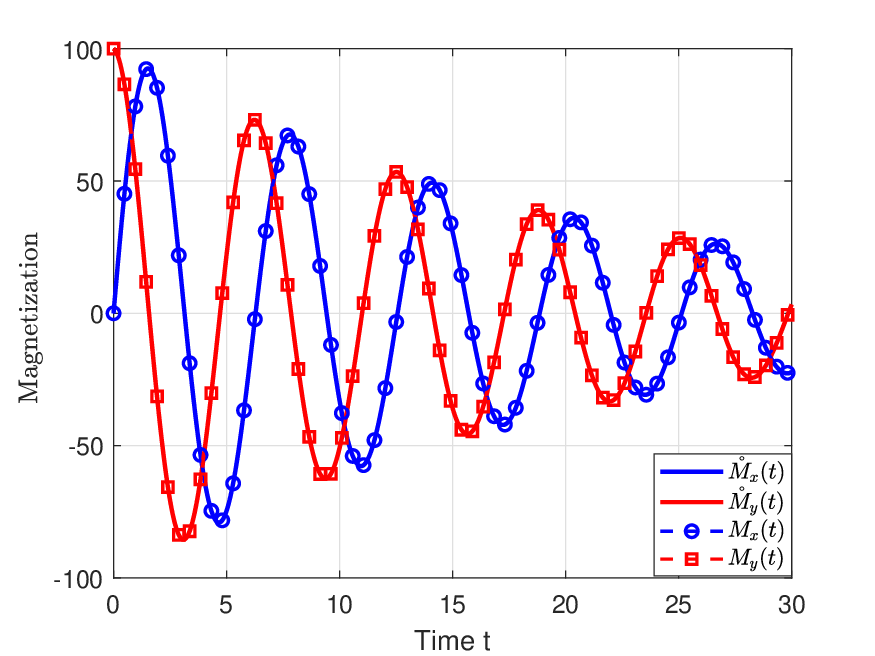}
    \caption{ $M_x(t)$ and $M_y(t)$ }
\end{subfigure}
\hfill
\begin{subfigure}[b]{0.48\textwidth}
    \centering
    \includegraphics[width=\linewidth]{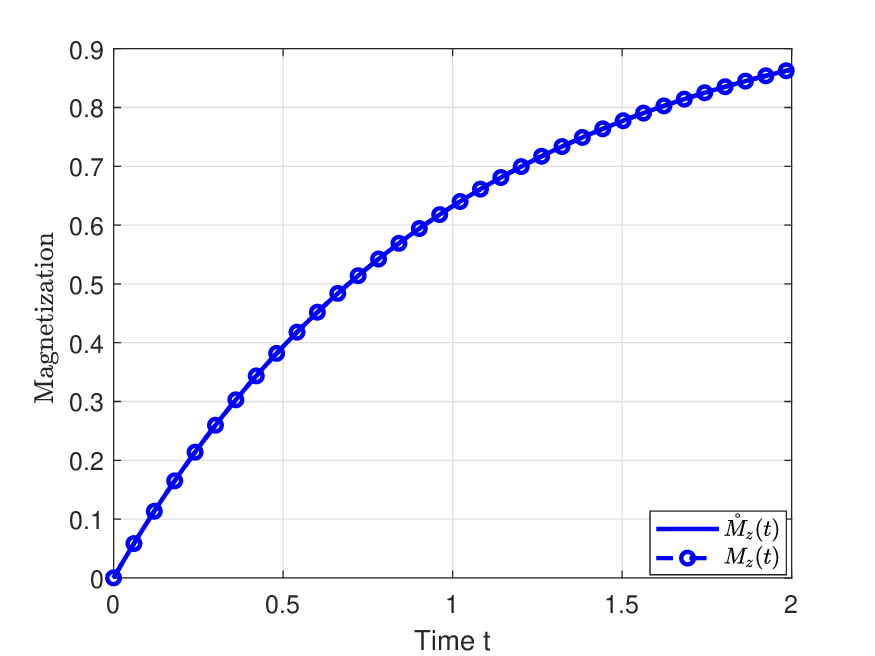}
    \caption{ $M_z(t)$}
\end{subfigure}

\caption{Comparison of exact and approximate solutions (with 50 terms in the series solution).}
\label{fig:bloch_comparison}

\end{figure}

Figure \ref{fig:bloch_comparison} demonstrate the comparison of our obtained solution of each magnetization components with the exact solution for the case $\alpha =1$. The obtained solution coincides exactly with the approximate solution on the entire time domain. Moreover, it is notable that our obtained solution shows the oscillatory damping behavior in the same extent as  of classical solution. This near-perfect overlap indicates that the proposed method accurately captures the transverse and longitudinal magnetic dynamics.

\begin{figure}[H]
    \centering
    \begin{subfigure}[H]{0.42\textwidth}
    \centering
    \includegraphics[width=\linewidth]{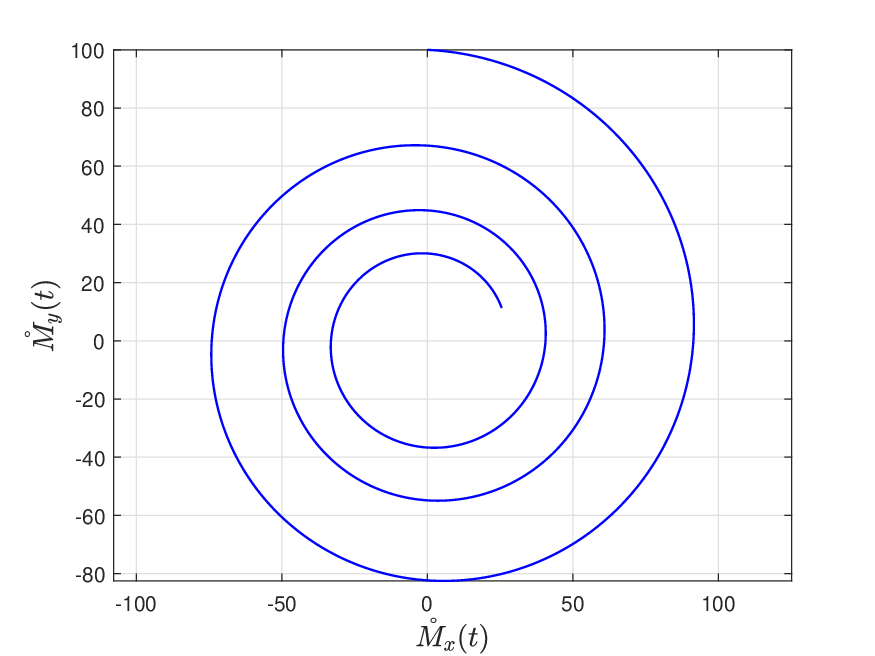}
    
    \caption{ $\mathring{M}_x(t)$ versus $\mathring{M}_y(t)$ at time $t=20$ and $\alpha =0.99$.} \label{M_xM_yapprox}
\end{subfigure}
\begin{subfigure}[H]{0.42\textwidth}
    \centering
    \includegraphics[width=\linewidth]{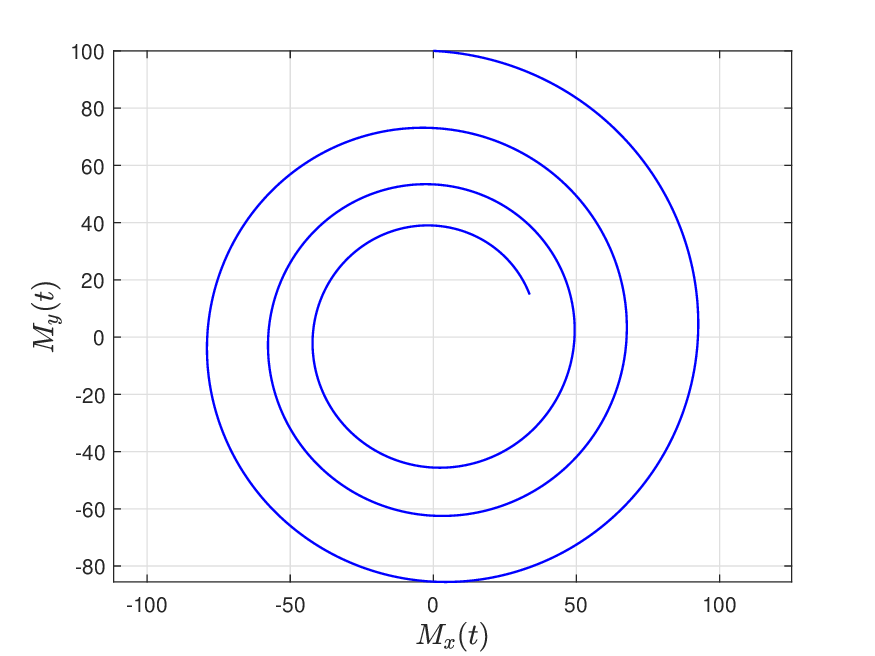}
    \caption{${M}_x(t)$ versus ${M}_y(t)$ at time $t=20$ and $\alpha =1$.}
    \label{M_xM_yexact}
    
\end{subfigure}
\\
    \begin{subfigure}[H]{0.42\textwidth}
        \centering
\includegraphics[width=\linewidth]{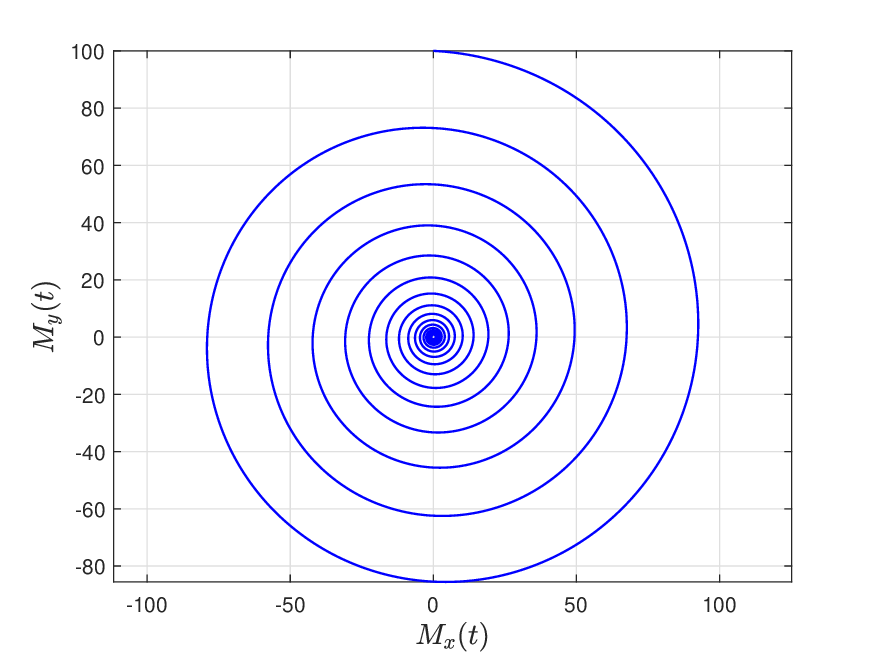}
        \caption{${M}_x(t)$ versus ${M}_y(t)$ at time $t=100$ and $\alpha =1$.}
        \label{M_xM_yexact_long}
    \end{subfigure}
\caption{ Phase plane representation of transverse components  (with 30 terms in the series solution).}
\label{fig:bloch_MxVsMy}
\end{figure}

\begin{figure}[H]
    \centering
    \begin{subfigure}{0.48\textwidth}
        \centering
\includegraphics[width=\linewidth]{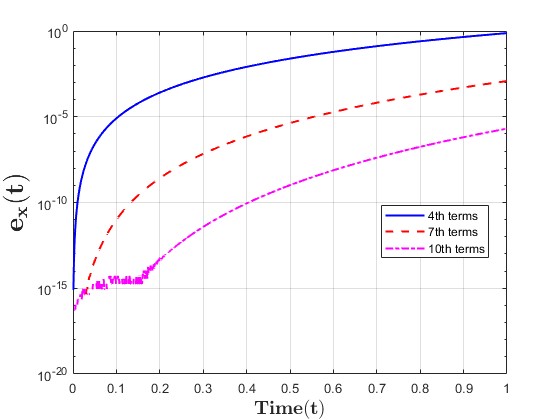}
        \caption{Absolute error $\mathbf{e}_x(t)$}
        \label{fig:error_mx}
    \end{subfigure}
    \hfill
    \begin{subfigure}{0.48   \textwidth}
        \centering
\includegraphics[width=\linewidth]{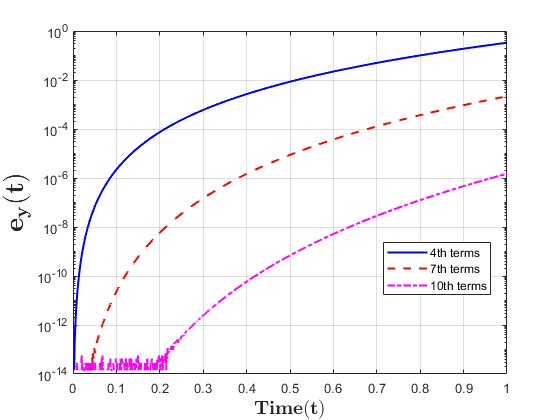}
        \caption{Absolute error $\mathbf{e}_y(t)$}
        \label{fig:error_my}
    \end{subfigure}
    \\
    \begin{subfigure}{0.5\textwidth}
        \centering
\includegraphics[width=\linewidth]{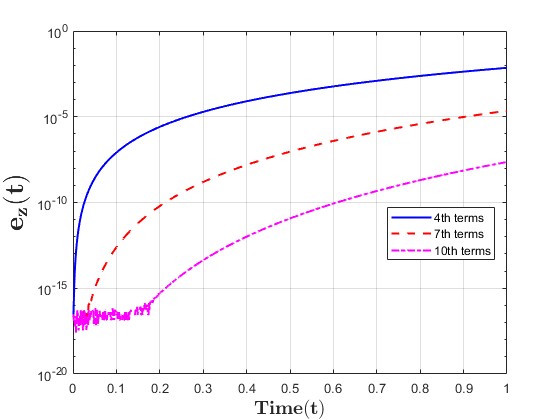}
        \caption{Absolute error $\mathbf{e}_z(t)$}
        \label{fig:error_mz}
    \end{subfigure}
    \caption{Comparison of absolute errors in $\mathring{M}_x(t)$, $\mathring{M}_y(t)$, and $\mathring{M}_z(t)$ with different number of terms in approximate solution.}
    \label{fig:absolute_errors}
\end{figure}

Figure \ref{fig:bloch_MxVsMy} illustrates the phase plane representation of transverse components. The phase plane trajectories of transverse components exhibits a clear spiraling pattern towards the origin. This indicate the damped oscillatory behavior governed by relaxation effects. The trajectories of obtained approximate solution is illustrated in Figure \ref{M_xM_yapprox}. The trajectories of analytical solution for $\alpha =1$ at different times $t=20$ and $t=100$ are demonstrated in Figure \ref{M_xM_yexact} and \ref{M_xM_yexact_long}, respectively. It is visible that the solution undergoes complete damping as time progress. Remarkably, we detect a slight deviation in the case of $\alpha =0.99$, where the decay is comparatively slower, highlighting the presence of memory effect from fractional formulation. Moreover, the smoothness of obtained trajectories along with their close agreement with the expected physical behavior validate the effectiveness of proposed approach.
We demonstrate the absolute error associated with each magnetization component in Figure \ref{fig:absolute_errors}. Here, we take approximate solutions up to the 4th, 7th, and 10th terms. The significant influence of the number of terms in the approximate solution on absolute errors is shown. We observe that the magnitude of error reduces remarkably as the number of terms increases, indicating convergence of the obtained solution. Moreover, the absolute error over the considered time interval is negligibly small, which highlighting the effectiveness and precision of the proposed method. 

\section{Conclusion} \label{sec6}
This work focuses on 
developing the approximate solution to the fractional Bloch equations using a hybrid approach, the Laplace-residual power series method.  This method integrates the Laplace transform with the conventional residual power series approach. The conversion of the fractional Bloch equation into a system of 
algebraic equations is carried out by operating the Laplace transform. In contrast to the conventional residual power series method, this series solution approach does not need the computation of fractional derivatives. Instead, this method depends on the simple limiting process to extract the unknown coefficients. Also, unlike the other approaches for solving the fractional Bloch equation, this method does not require perturbation or discretizations. The obtained solution keeps a good agreement with the exact solution for all magnetization components. The absolute and relative error analyses demonstrate the high level of accuracy of the suggested method. Furthermore, the phase-plane trajectory illustrates the evolution of the magnetization vector and highlight the impact of fractional order parameter on the relaxation dynamics of the Bloch system. These results indicate that LRPSM is an efficient and reliable technique for studying fractional Bloch models encountered in magnetic resonance phenomena. 

\section*{Acknowledgments}
The authors would like to sincerely thank the anonymous editors and reviewers for their suggestions. Varsha R is supported by the University Grant Commission of India.

\bibliographystyle{unsrt}
\bibliography{Ref}

\end{document}